# Intuitive Considerations Clarifying the Origin and Applicability of the Benford Law

G. Whyman[*], E. Shulzinger, Ed. Bormashenko

Ariel University, Faculty of Natural Sciences, Department of Physics, Ariel, P.O.B.3, 40700, Israel

## Abstract

The diverse applications of the Benford law attract investigators working in various fields of physics, biology and sociology. At the same time, the groundings of the Benford law remain obscure. Our paper demonstrates that the Benford law arises from the positional (place-value) notation accepted for representing various sets of data. An alternative to Benford formulae to predict the distribution of digits in statistical data are derived. Application of these formulae to the statistical analysis of infrared spectra of polymers is presented. Violations of the Benford Law are discussed.

## Introduction

The Benford law is a phenomenological, contra-intuitive law observed in many naturally occurring tables of numerical data; also called the first-digit law, first digit phenomenon, or leading digit phenomenon. It states that in listings, tables of statistics, etc., the digit 1 tends to occur with a probability of 30%, much greater than the expected value of 11.1% (i.e. one digit out of 9) [1-2]. The discovery of Benford's law goes back to 1881, when the American astronomer Simon Newcomb noticed that in logarithm tables (used at that time to perform calculations), the earlier pages (which contained numbers that started with 1) were much more worn and smudged than the later pages. Newcomb noted, "that the ten digits do not occur with equal frequency must be evident to any making use of logarithmic tables, and noticing how much faster first pages wear out than the last ones [1]." The phenomenon was re-discovered by the physicist Frank Benford, who tested it on data extracted from 20 different domains, as different as surface areas of rivers, physical constants, molecular weights, etc. Since then, the law has been credited to Benford [2]. The Benford law is expressed by the following statement: the occurrence of first significant digits in data sets follows a logarithmic distribution:

$$P(n) = \log_{10}\left(1 + \frac{1}{n}\right), \qquad n = 1, 2, ..., 9 \qquad (1)$$

where $P(n)$ is the probability of a number having the first non-zero digit $n$.

Since its formulation, Benford's law has been applied for the analysis of a broad variety of statistical data, including atomic spectra [3], population dynamics [4], magnitude and depth of earthquakes [5], genomic data [6-7], mantissa distributions of pulsars [8], and economic data [9-10]. While Benford's law definitely applies to many situations in the real world, a satisfactory explanation has been given only recently through the works of Hill et al. [11-13], who called the Benford distribution "the law of statistical folklore". Important intuitive physical insights in the grounding of the Benford law, relating its origin to the scaling invariance of physical laws, were reported by Pietronero et al. [14]. Engel et al. demonstrated that the Benford law takes place approximatively for exponentially distributed numbers [15]. Fewster supplied a simple "geometrical" reasoning of the Benford law [16]. The breakdown of the Benford law was reported for certain sets of statistical data [17-19].

It should be mentioned that the grounding and applicability of the Benford law remain highly debatable [13]. In spite of this, the Benford law was effectively exploited for detecting fraud in accounting data [18]. Quantifying non-stationarity effects on organization of turbulent motion by Benford's law was reported recently [20]. Our paper supplies intuitive reasoning clarifying the origin of the Benford law.

1. **New Results**
   1.1. **The Origin of the Benford Law, and the Positional (place-value) notation**

In practice, measured quantities or analyzed data are restricted by a prescribed accuracy defined by a number of significant digits. This means that mantissas of decimal numbers, which are simply integers, are restricted from above by some integer, say, $m+1$.

Taking in mind the above mentioned, consider a set $\{1,2,...,m\}$. When $m \to \infty$, this set coincides with the full set of integers. Let us elucidate how the frequency $f_n(m)$ of numbers beginning with the digit 1 ($n=1$) depends on $m$. In the first 6 lines of Table 1, the examples for the values of $m$ are presented for which $f_1(m)$ successively reaches minimum and maximum. It is seen that the above frequency changes quasi-periodically with increasing $m$, decreasing and increasing, and reaches its minima and maxima in turn for selected values of $m$ (see Figure 1). Successive minimums $f_{\min,n}(k)$ and maximums $f_{\max,n}(k)$ are enumerated by $k=1,2....$

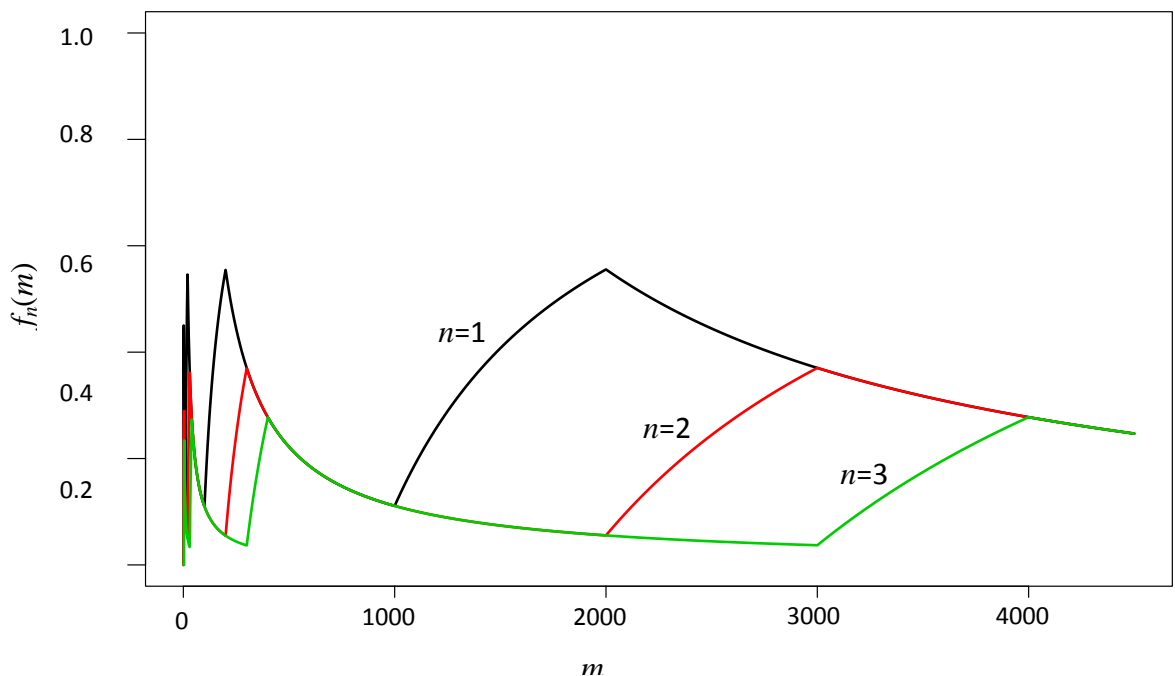

Figure 1. The dependence of first-digit frequency on the upper mantissa limit for the digits $n=1$ (black), 2 (red), and 3 (green).

As another example, in the following lines of Table 1, the minimal and maximal frequencies $f_{\min,5}(k)$, $f_{\max,5}(k)$ and $f_{\min,9}(k)$, $f_{\max,9}(k)$ of integers beginning with the digits 5 and 9 are given. It is seen that the maximal and minimal frequencies decrease for the sequence $n=1,5,9$: the number of integers beginning with these digits remains the same, but the sizes of the corresponding intervals $[1,m]$ grow (compare $m$ in the third column for different $n$ and the same $k$).

Table 1. Frequencies of integers beginning with different figures.

| First digit, $n$, of the number | $m$ | $\{1,2,...,m\}$ | Amount, $p$, of numbers | $k$ | Minimal and maximal frequencies, $p/m$ |
|---|---|---|---|---|---|
| 1 | 9 | 1,2,…,9 | 1 | 1 | 1/9 |
| | 19 | 1,2,..,19 | 11 | | 11/19 |
| | 99 | 1,2,…,99 | 11 | 2 | 11/99=1/9 |
| | 199 | 1,2,…,199 | 111 | | 111/199 |
| | 999 | 1,2,…,999 | 111 | 3 | 111/999=1/9 |
| | 1999 | 1,2,…,1999 | 1111 | | 1111/1999 |
| 5 | 49 | 1,2,…,49 | 1 | 1 | 1/49 |
| | 59 | 1,2,…,59 | 11 | | 11/59 |
| | 499 | 1,2,…,499 | 11 | 2 | 11/499 |
| | 599 | 1,2,…,599 | 111 | | 111/599 |
| | 4999 | 1,2,…,4999 | 111 | 3 | 111/4999 |
| | 5999 | 1,2,…,5999 | 1111 | | 1111/5999 |
| 9 | 89 | 1,2,…,89 | 1 | 1 | 1/89 |
| | 99 | 1,2,…,99 | 11 | | 11/99 |
| | 899 | 1,2,…,899 | 11 | 2 | 11/899 |
| | 999 | 1,2,…,999 | 111 | | 111/999 |
| | 8999 | 1,2,…,8999 | 111 | 3 | 111/8999 |
| | 9999 | 1,2,…,9999 | 1111 | | 1111/9999 |

As is seen from Table 1, the successive minima and maxima, enumerated by $k$, may be written as

$$f_{\min,n}(k) = \frac{1 \cdot 10^{k-1} + 1 \cdot 10^{k-2} + \cdots + 1}{(n-1) \cdot 10^k + 9 \cdot 10^{k-1} + 9 \cdot 10^{k-2} + \cdots + 9}, \quad (2)$$

$$f_{\max,n}(k) = \frac{1 \cdot 10^k + 1 \cdot 10^{k-1} + \cdots + 1}{n \cdot 10^k + 9 \cdot 10^{k-1} + 9 \cdot 10^{k-2} + \cdots + 9} \quad (3)$$

for $k = 1,2,3\ldots$. All the sums in (2), (3) are calculated as sums of the geometric sequence

$$f_{\min,n}(k) = \frac{10^k - 1}{9(n \cdot 10^k - 1)}, \quad (4)$$

$$f_{\max,n}(k) = \frac{10^{k+1} - 1}{9[(n+1) \cdot 10^k - 1]}. \quad (5)$$

Letting $k$ go to infinity (which also means letting corresponding values of $m$ in Table 1 to go to infinity), results in

$$f_{\min,n} = \lim_{k \to \infty} f_{\min,n}(k) = \frac{1}{9n}, \quad (6)$$

$$f_{\max,n} = \lim_{k \to \infty} f_{\max,n}(k) = \frac{10}{9(n+1)}. \quad (7)$$

The probability of the occasional choosing of a particular number beginning with the digit $n$ from the whole set of integers may be estimated as a normalized arithmetic mean or a normalized geometric mean of the minimal (6) and maximal (7) frequencies:

$$P_{\text{arith}}(n) = [f_{\min,n} + f_{\max,n}] / \left( \sum_{i=1}^{9} (f_{\min,i} + f_{\max,i}) \right)$$

$$P_{\text{geom}}(n) = \sqrt{f_{\min,n} f_{\max,n}} / \left( \sum_{i=1}^{9} \sqrt{f_{\min,i} \cdot f_{\max,i}} \right).$$

The final result is

$$P_{\text{arith}}(n) = \frac{\frac{10}{n+1} + \frac{1}{n}}{\sum_{i=1}^{9} \left( \frac{10}{i+1} + \frac{1}{i} \right)} \tag{8}$$

$$P_{\text{geom}}(n) = \frac{1}{\sqrt{n(n+1)} \sum_{i=1}^{9} 1/\sqrt{i(i+1)}} . \tag{9}$$

The results of equations (8) and (9) are compared with the Benford formula (1) in Table 2 and Figure 2. As is seen, the normalized geometric mean shows very good agreement, even though the mathematical forms of (1) and (9) are different.

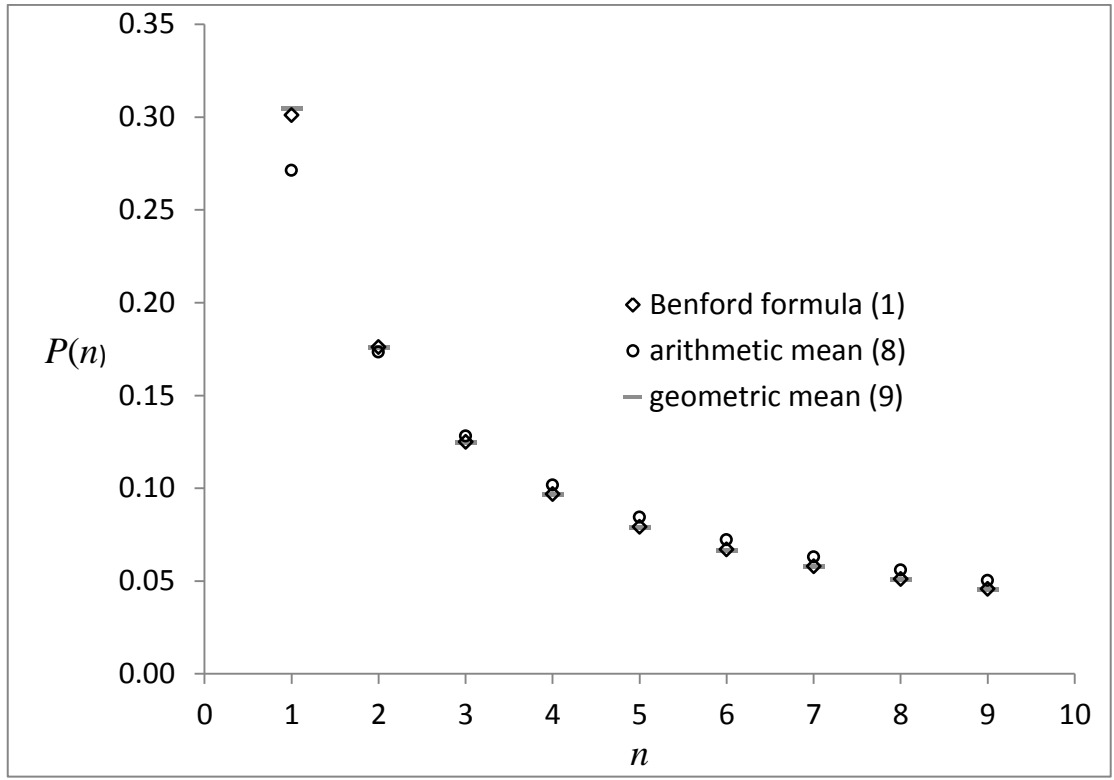

Figure 2. Comparison of equations (8) and (9) with the Benford formula (1).

Table 2. Comparison of equations (8) and (9) with the Benford formula (1).

| $n$ | 1 | 2 | 3 | 4 | 5 | 6 | 7 | 8 | 9 |
|---|---|---|---|---|---|---|---|---|---|
| Benford | 0.3010 | 0.1761 | 0.1249 | 0.09691 | 0.07918 | 0.06695 | 0.05799 | 0.05115 | 0.04576 |
| Geometric mean, Eq. (9) | 0.3046 | 0.1759 | 0.1244 | 0.09632 | 0.07865 | 0.06647 | 0.05756 | 0.05077 | 0.04541 |
| Arithmetic mean, Eq. (8) | 0.2712 | 0.1733 | 0.1281 | 0.1017 | 0.08439 | 0.07212 | 0.06297 | 0.05589 | 0.05023 |

The results (6)-(9) allow an obvious generalization for the case of an arbitrary base $N$ of the positional digit system:

$$f^N_{\min,n} = \frac{1}{(N-1)n}, \qquad f^N_{\max,n} = \frac{N}{(N-1)(n+1)}$$

$$P^N_{\text{arith}}(n) = \left(\frac{N}{n+1} + \frac{1}{n}\right) / \sum_{i=1}^{N-1}\left(\frac{N}{i+1} + \frac{1}{i}\right) \tag{10}$$

$$P^N_{\text{geom}}(n) = \frac{1}{\sqrt{n(n+1)}\,\sum_{i=1}^{N-1} 1/\sqrt{i(i+1)}} \tag{11}$$

where $1 \leq n \leq N-1$. In particular, in the binary system ($N=2$), all the right-hand sides of the four last equations turn to 1 for $n=1$ (all the numbers presented in the binary system begin with 1).

It is well known that in many cases the Benford distribution does not hold. This may happen, e.g., under some restriction on the set of admissible numbers. For example, if the inequality $1 \leq l < 1000$ is imposed on the random sample of integers $l$ (or mantissas of real numbers), the probability $P(1)$ will be close to 1/9 (see Table 1), and not to the value predicted by the Benford formula or by equations (8, 9), which is about 3 times larger. More generally, the necessary condition is that the set $\{1, 2, \ldots, m\}$ to which a random sample of integers belong should contain the same numbers of minimal (4) and maximal (5) frequencies. In any case, if some restrictions take place, the following inequalities should be fulfilled:

$$f_{\min,n} \leq P(n) \leq f_{\max,n}$$

or

$$\frac{1}{9n} \leq P(n) \leq \frac{10}{9(n+1)} \tag{12}$$

in the decimal system. In digit systems with a lower base $N$, the appropriate inequalities are stronger:

$$\frac{1}{(N-1)n} \leq P(n) \leq \frac{N}{(N-1)(n+1)}.$$

A favorable situation for the Benford distribution appears when admissible numbers belong to a function range in a vicinity of infinite singularity. In this case, restrictions on *m* are absent, and the statement of tending *m* to infinity in (6) and (7) becomes reasonable.

### 1.2. Exemplification of New Results: Applicability of the Obtained Results to the Analysis of Infrared Spectra of Polymers

In our recent paper we demonstrated that the Benford law takes place within the absorbance domain of infrared (IR) spectra of polymers [21]. The IR spectra may be treated as sets of values of absorbance corresponding to the sets of wavenumbers. Consider now validity of Eqs. 8-9 to the actual distribution of leading digits in the absorbance spectra of polymers studied in Ref. 21, and represented in Fig. 3.It is recognized that the geometrical averaging given by Eq. 11 supplies the best correspondence with the experimental results.

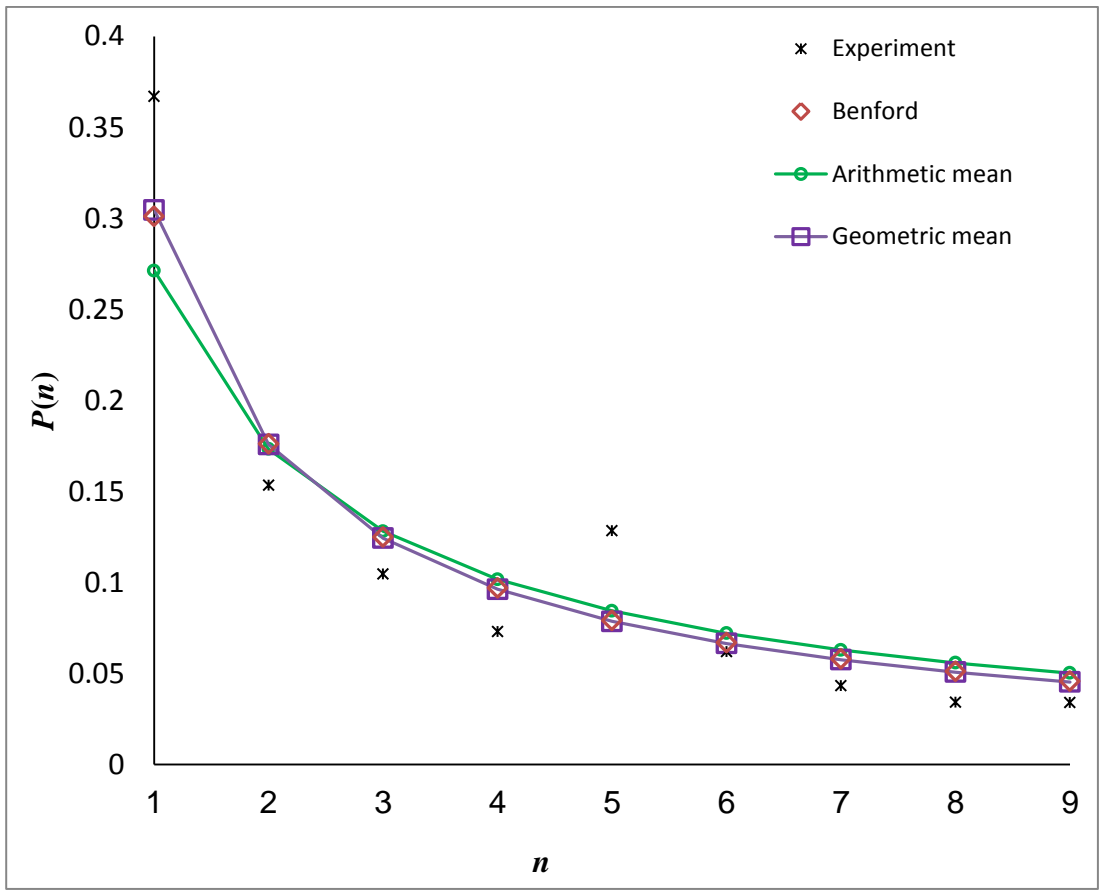

Figure 3. The actual frequencies of leading digits appearing in the set of absorbance spectra vs. the Benford law and Eqs. 10 and 11. The correlation coefficients are: $R=0.964$ for the Benford distribution, $R=0.956$ for Eq. 10 (arithmetic average approximation) and $R=0.966$ for Eq. 11 (geometrical average approximation).

**Summary**

The present article places emphasis on the Benford law as a consequence of the structure of positional digit systems. From this point of view, attempts of explanations based on scale invariance, base invariance or even representing of the Benford law as a mysterious law of nature, at least call for refinement. A very convincing example is the binary positional digit system (with a base of 2) for which the "Benford law" should state that the probability of finding the digit 1 at the first place of a number is 100%.

As shown above, some statistical estimation of the probability of finding the digits at the first place of a number can be given, which obeys a different mathematical form alternative to the Benford law. This form, which is expressed by the derived equations (10), (11), gives practically the same numerical results as the Benford formula.

Limitations from below and from above on admissible numbers imposed a priori lead to violations of the Benford law. Some inequalities concerning these violations can be useful.

**Acknowledgements**

GW thanks to Israel Ministry of Absorption for years-long generous support and to his sister Elena Vaiman for her help.